\documentclass[12pt,reqno,a4wide]{amsart}

\allowdisplaybreaks

\newtheorem{theorem}{Theorem}

\newtheorem{remark}[theorem]{Remark}

\begin{document}

	\title[On the generating function of $p$-Bernoulli numbers]{On the generating function of $p$-Bernoulli numbers: an alternative approach}

	\author[L. Karg{\i}n]{Levent Karg{\i}n}
	\address{(Karg{\i}n) Vocational School of Technical Sciences, Akdeniz University, Antalya TR-07058 Turkey}
	\email{leventkargin48@gmail.com}

	\author[M. Rahmani]{Mourad Rahmani}
	\address{(Rahmani) Faculty of Mathematics, USTHB, PO. Box 32, El Alia, 16111\\
                    Algiers, Algeria}
	\email{rahmani.mourad@gmail.com}
	
	
	\begin{abstract}
		In this note, we give an alternative proof of the generating function of $p$-Bernoulli numbers. Our argument is based on the Euler's integral representation.
	\end{abstract}
	

	\maketitle

	\section{Introduction}
The $p$-Bernoulli numbers $B_{n, p}\ (n,p \geq 0)$ was introduced by the second author in a recent paper \cite{Rahmani} as a generalization of Bernoulli numbers, which are defined by the following exponential generating function
\begin{equation}
\sum_{n\geq 0}B_{n, p}\frac{t^{n}}{n!}=\ _{2}F_{1}\left(  1,1;p+2;1-e^{t}\right),
\end{equation}
where $_{2}F_{1}\left(  a,b;c;z\right)  $ denotes the Gaussian hypergeometric function. Please refer to \cite{Kargin,Kargin1,doi:10.2989/16073606.2017.1418762,Rahmani} for more details on these numbers.

In \cite{doi:10.2989/16073606.2017.1418762} the authors gave the closed formula for the exponential generating function of $B_{n, p}$ in terms of the harmonic numbers.
\begin{theorem}
\label{thm}
For $p\geq 0$
\begin{equation}
\label{eq:main}
\sum_{n=0}^{\infty}B_{n, p}\frac{t^{n}}{n!} = \frac{(p+1)(t-H_{p})e^{pt}}{(e^{t}-1)^{p+1}}+(p+1)\sum_{k=1}^{p}\binom{p}{k}\frac{H_{k}}{(e^{t}-1)^{k+1}},
\end{equation}
where $H_{n}$ is the harmonic numbers defined by
\begin{equation*}
H_{0} =0,\  H_{n} = \sum_{j=1}^{n}\frac{1}{j}\ \ (n\geq 1).
\end{equation*}
\end{theorem}

In this note, we give an alternative proof of the generating function of $p$-Bernoulli numbers. The proof of Theorem \ref{thm}, is based upon the formula \cite[p.361]{Rahmani}
\begin{equation}%
{\displaystyle\sum\limits_{n\geq0}}
B_{n,p}\frac{t^{n}}{n!}   =\left(  p+1\right)
{\displaystyle\int\limits_{0}^{1}}
\frac{\left(  1-x\right)  ^{p}}{1-\left(  1-e^{t}\right)  x}dx. \label{ZAE}
\end{equation}
\begin{remark}
Note that Prodinger and Selkirk in \cite{Prodinger} presented a simple proof of (\ref{eq:main}) by using elementary methods. Another proof was also given by Kuba in \cite{Kuba}.
\end{remark}
\section{Proof}
From (\ref{ZAE}) and binomial formula, we get%
\[%
{\displaystyle\sum\limits_{n\geq0}}
B_{n,p}\frac{t^{n}}{n!}=\left(  p+1\right)
{\displaystyle\sum\limits_{k=0}^{p}}
\dbinom{p}{k}%
{\displaystyle\int\limits_{0}^{1}}
\frac{\left(  -x\right)  ^{k}}{1-\left(  1-e^{t}\right)  x}dx.
\]
Now, making the substitution $u=\left(  1-e^{t}\right)  x$, we have%
\begin{align*}%
{\displaystyle\sum\limits_{n\geq0}}
B_{n,p}\frac{t^{n}}{n!} &  =\left(  p+1\right)
{\displaystyle\sum\limits_{k=0}^{p}}
\dbinom{p}{k}\frac{1}{\left(  e^{t}-1\right)  ^{k+1}}%
{\displaystyle\int\limits_{0}^{1-e^{t}}}
\frac{-u^{k}}{1-u}du\\
&  =\left(  p+1\right)
{\displaystyle\sum\limits_{k=0}^{p}}
\dbinom{p}{k}\frac{1}{\left(  e^{t}-1\right)  ^{k+1}}\left(
{\displaystyle\int\limits_{0}^{1-e^{t}}}
\frac{1-u^{k}}{1-u}du-%
{\displaystyle\int\limits_{0}^{1-e^{t}}}
\frac{1}{1-u}du\right)  \\
&  =\frac{\left(  p+1\right)  te^{pt}}{\left(  e^{t}-1\right)  ^{p+1}}+\left(
p+1\right)
{\displaystyle\sum\limits_{k=0}^{p}}
\dbinom{p}{k}\frac{1}{\left(  e^{t}-1\right)  ^{k+1}}%
{\displaystyle\int\limits_{0}^{1-e^{t}}}
\frac{1-u^{k}}{1-u}du.
\end{align*}
Since the integral representation of harmonic numbers is
\[
H_{n}=%
{\displaystyle\int\limits_{0}^{1}}
\frac{1-x^{n}}{1-x}dx.
\]
Then we have%
\begin{align*}%
{\displaystyle\sum\limits_{n\geq0}}
B_{n,p}\frac{t^{n}}{n!} &  =\frac{\left(  p+1\right)  te^{pt}}{\left(
e^{t}-1\right)  ^{p+1}}+\left(  p+1\right)
{\displaystyle\sum\limits_{k=0}^{p}}
\dbinom{p}{k}\frac{1}{\left(  e^{t}-1\right)  ^{k+1}}\left(  H_{k}-%
{\displaystyle\int\limits_{1-e^{t}}^{1}}
\frac{1-u^{k}}{1-u}du\right)  \\
&  =\frac{\left(  p+1\right)  te^{pt}}{\left(  e^{t}-1\right)  ^{p+1}}+\left(
p+1\right)
{\displaystyle\sum\limits_{k=0}^{p}}
\dbinom{p}{k}\frac{1}{\left(  e^{t}-1\right)  ^{k+1}}\left(  H_{k}-%
{\displaystyle\sum\limits_{i=0}^{k-1}}
{\displaystyle\int\limits_{1-e^{t}}^{1}}
u^{i}du\right)  \\
&  =\frac{\left(  p+1\right)  te^{pt}}{\left(  e^{t}-1\right)  ^{p+1}}+\left(
p+1\right)
{\displaystyle\sum\limits_{k=0}^{p}}
\dbinom{p}{k}\frac{1}{\left(  e^{t}-1\right)  ^{k+1}}\left(  H_{k}-%
{\displaystyle\sum\limits_{i=1}^{k}}
\frac{1-\left(  1-e^{t}\right)  ^{i}}{i}\right).
\end{align*}
Using again the binomial formula, we obtain%
\begin{multline}%
{\displaystyle\sum\limits_{n\geq0}}
B_{n,p}\frac{t^{n}}{n!}=\frac{\left(  p+1\right)  te^{pt}}{\left(
e^{z}-1\right)  ^{p+1}}+\left(  p+1\right)
{\displaystyle\sum\limits_{k=0}^{p}}
\dbinom{p}{k}\frac{H_{k}}{\left(  e^{t}-1\right)  ^{k+1}}\label{LLM1}\\
-\frac{\left(  p+1\right)  e^{pt}}{\left(  e^{z}-1\right)  ^{p+1}}%
{\displaystyle\sum\limits_{k=0}^{p}}
\dbinom{p}{k}\frac{\left(  e^{t}-1\right)  ^{p-k}}{e^{pt}}%
{\displaystyle\sum\limits_{i=1}^{k}}
\frac{1}{i}\left(  1-%
{\displaystyle\sum\limits_{s=0}^{i}}
\dbinom{i}{s}\left(  -1\right)  ^{s}e^{st}\right).
\end{multline}
Therefore%
\begin{multline}
{\sum\limits_{n\geq0}}B_{n,p}\frac{t^{n}}{n!}=\frac{\left(  p+1\right)
te^{pt}}{\left(  e^{t}-1\right)  ^{p+1}}+\left(  p+1\right)  {\sum
\limits_{k=0}^{p}}\dbinom{p}{k}\frac{H_{k}}{\left(  e^{t}-1\right)  ^{k+1}%
}\label{aqwer}\\
+\frac{\left(  p+1\right)  e^{pt}}{\left(  e^{t}-1\right)  ^{p+1}}%
{\sum\limits_{k=0}^{p}}\dbinom{p}{k}\left(  e^{t}-1\right)  ^{p-k}%
{\sum\limits_{s=1}^{i}}\left(  -1\right)  ^{s}e^{\left(  s-p\right)  t}%
{\sum\limits_{i=1}^{k}}\frac{1}{i}\dbinom{i}{s}.
\end{multline}
Since%
\[%
%TCIMACRO{\dsum \limits_{i=1}^{k}}%
%BeginExpansion
{\displaystyle\sum\limits_{i=1}^{k}}
%EndExpansion
\frac{1}{i}\dbinom{i}{s}=\frac{1}{s}\dbinom{k}{s},
\]
then (\ref{aqwer}) becomes%
\begin{multline}%
%TCIMACRO{\dsum \limits_{n\geq0}}%
%BeginExpansion
{\displaystyle\sum\limits_{n\geq0}}
%EndExpansion
B_{n,p}\frac{t^{n}}{n!}=\frac{\left(  p+1\right)  te^{pt}}{\left(
e^{t}-1\right)  ^{p+1}}+\left(  p+1\right)
%TCIMACRO{\dsum \limits_{k=0}^{p}}%
%BeginExpansion
{\displaystyle\sum\limits_{k=0}^{p}}
%EndExpansion
\dbinom{p}{k}\frac{H_{k}}{\left(  e^{t}-1\right)  ^{k+1}}\label{LL1}\\
+\frac{\left(  p+1\right)  e^{pt}}{\left(  e^{t}-1\right)  ^{p+1}}%
%TCIMACRO{\dsum \limits_{s=1}^{p}}%
%BeginExpansion
{\displaystyle\sum\limits_{s=1}^{p}}
%EndExpansion
\text{ }\left(  -1\right)  ^{s}\frac{1}{s}e^{\left(  s-p\right)  t}\left(
%TCIMACRO{\dsum \limits_{k=0}^{p}}%
%BeginExpansion
{\displaystyle\sum\limits_{k=0}^{p}}
%EndExpansion
\dbinom{k}{s}\dbinom{p}{k}\left(  e^{t}-1\right)  ^{p-k}\right).
\end{multline}
It is easily verified that
\[%
%TCIMACRO{\dsum \limits_{k=0}^{p}}%
%BeginExpansion
{\displaystyle\sum\limits_{k=0}^{p}}
%EndExpansion
\dbinom{p}{k}\dbinom{k}{s}\left(  e^{t}-1\right)  ^{p-k}=\dbinom{p}%
{s}e^{\left(  p-s\right)  t}.
\]
Thus we have%
\begin{equation}%
%TCIMACRO{\dsum \limits_{n\geq0}}%
%BeginExpansion
{\displaystyle\sum\limits_{n\geq0}}
%EndExpansion
B_{n,p}\frac{t^{n}}{n!}=\frac{\left(  p+1\right)  e^{pt}}{\left(
e^{t}-1\right)  ^{p+1}}\left(  t+%
%TCIMACRO{\dsum \limits_{s=1}^{p}}%
%BeginExpansion
{\displaystyle\sum\limits_{s=1}^{p}}
%EndExpansion
\text{ }\left(  -1\right)  ^{s}\frac{1}{s}\dbinom{p}{s}\right)  +\left(
p+1\right)
%TCIMACRO{\dsum \limits_{k=0}^{p}}%
%BeginExpansion
{\displaystyle\sum\limits_{k=0}^{p}}
%EndExpansion
\dbinom{p}{k}\frac{H_{k}}{\left(  e^{t}-1\right)  ^{k+1}}.\label{LLZ}%
\end{equation}
Formula (\ref{eq:main}) now follows from (\ref{LLZ}) and
\[
H_{p}=-%
%TCIMACRO{\dsum \limits_{s=1}^{p}}%
%BeginExpansion
{\displaystyle\sum\limits_{s=1}^{p}}
%EndExpansion
\text{ }\left(  -1\right)  ^{s}\frac{1}{s}\dbinom{p}{s}.
\]
The proof of Theorem \ref{thm} is complete.
\begin{remark} 
The polynomial version of the generating function of $p$-Bernoulli numbers is easily obtained by multiplying the right-hand side of (\ref{eq:main}) by $e^{xt}$.
\[
\label{eq:main1}
\sum_{n=0}^{\infty}B_{n, p}(x)\frac{t^{n}}{n!} = \frac{(p+1)(t-H_{p})e^{(x+p)t}}{(e^{t}-1)^{p+1}}+(p+1)e^{xt}\sum_{k=1}^{p}\binom{p}{k}\frac{H_{k}}{(e^{t}-1)^{k+1}}.
\]
\end{remark}

\bibliographystyle{plain}

\end{document}